\documentclass{amsart}
\usepackage{graphicx}
\usepackage{amscd}
\usepackage{amsmath}
\usepackage{amsfonts}
\usepackage{amssymb}
\newtheorem{theorem}{Theorem}
\theoremstyle{plain}

\newtheorem{corollary}{Corollary}

\newtheorem{remark}{Remark}

\numberwithin{equation}{section}

\begin{document}
\title[Logarithmic Sobolev inequalities]{Isoperimetry and Symmetrization for Logarithmic Sobolev inequalities}
\author{Joaquim Mart\'{\i}n$^{\ast}$}
\address{Department of Mathematics\\
Universitat Aut\`{o}noma de Barcelona}
\email{jmartin@mat.uab.cat}
\author{Mario Milman}
\address{Department of Mathematics\\
Florida Atlantic University}
\email{extrapol@bellsouth.net}
\urladdr{http://www.math.fau.edu/milman}
\thanks{$^{\ast}$ Supported in part by MTM2004-02299 and by CURE 2005SGR00556}
\thanks{This paper is in final form and no version of it will be submitted for
publication elsewhere.}

\begin{abstract}
Using isoperimetry and symmetrization we provide a unified framework to study
the classical and logarithmic Sobolev inequalities. In particular, we obtain
new Gaussian symmetrization inequalities and connect them with logarithmic
Sobolev inequalities. Our methods are very general and can be easily adapted
to more general contexts.
\end{abstract}\maketitle
\tableofcontents

\section{Introduction}

The classical $L^{2}$-Sobolev inequality states that%
\[
\left|  \nabla f\right|  \in L^{2}(\mathbb{R}^{n})\Rightarrow f\in
L^{p_{n}^{\ast}}(\mathbb{R}^{n}),\text{ where }\frac{1}{p_{n}^{\ast}}=\frac
{1}{2}-\frac{1}{n}.
\]
Consequently, $\lim_{n\rightarrow\infty}p_{n}^{\ast}=2$ and, therefore, the
improvement on the integrability of $f$ disappears as $n\rightarrow\infty$. On
the other hand, Gross \cite{gr} showed that, if one replaces $dx$ by the
Gaussian measure $d\gamma_{n}(x)=(2\pi)^{-n/2}e^{-|x|^{2}/2}dx,$ we have%

\begin{equation}
\int\left|  f(x)\right|  ^{2}\ln\left|  f(x)\right|  d\gamma_{n}(x)\leq
\int\left|  \nabla f(x)\right|  ^{2}d\gamma_{n}(x)+\left\|  f\right\|
_{2}^{2}\ln\left\|  f\right\|  _{2}. \label{launo}%
\end{equation}
This is Gross' celebrated logarithmic Sobolev inequality (= lS inequality),
the starting point of a new field, with many important applications to PDEs,
Functional Analysis, Probability, etc. (as a sample, and only a sample, we
mention \cite{be}, \cite{dav}, \cite{le1}, \cite{log}, and the references
therein). The inequality (\ref{launo}) gives a logarithmic improvement on the
integrability of $f$, with constants independent of $n,$ that persists as
$n\mapsto\infty,$ and is best possible$.$ Moreover, rescaling (\ref{launo})
leads to $L^{p}$ variants of this inequality, again with constants independent
of the dimension (cf. \cite{gr}),
\[
\int\left|  f(x)\right|  ^{p}\ln\left|  f(x)\right|  d\gamma_{n}(x)\leq
\frac{p}{2(p-1)}\operatorname{Re}<Nf,f_{p}>+\left\|  f\right\|  _{p}^{p}%
\ln\left\|  f\right\|  _{p},
\]
where $<f,g>$ $=$ $\int f\bar{g}d\gamma_{n},$ $<Nf,f>$ $=\int\left|  \nabla
f(x)\right|  ^{2}d\gamma_{n}(x),$ $f_{p}=(sgn(f))\left|  f\right|  ^{p-1}.$

In a somewhat different direction, Feissner's thesis \cite{fe} under Gross,
takes up the embedding implied by (\ref{launo}), namely
\[
W_{2}^{1}(\mathbb{R}^{n},d\gamma_{n})\subset L^{2}(LogL)(\mathbb{R}%
^{n},d\gamma_{n}),
\]
where the norm of $W_{2}^{1}(\mathbb{R}^{n},d\gamma_{n})$ is given by
\[
\left\|  f\right\|  _{W_{2}^{1}(\mathbb{R}^{n},d\gamma_{n})}=\left\|  \nabla
f\right\|  _{L^{2}(\mathbb{R}^{n},d\gamma_{n})}+\left\|  f\right\|
_{L^{2}(\mathbb{R}^{n},d\gamma_{n})},
\]
and extends it to $L^{p},$ even Orlicz spaces. A typical result\footnote{For
the most part the classical work on functional lS inequalities has focussed on
$L^{2},$ or more generally, $L^{p}$ and Orlicz spaces.} from \cite{fe} is
given by
\begin{equation}
W_{p}^{1}(\mathbb{R}^{n},d\gamma_{n})\subset L^{p}(LogL)^{p/2}(\mathbb{R}%
^{n},d\gamma_{n}),\; 1<p<\infty. \label{launodos}%
\end{equation}

The connection between lS inequalities and the classical Sobolev estimates has
been investigated intensively. For example, it is known that (\ref{launo})
follows from the classical Sobolev estimates with sharp constants (cf.
\cite{be1}, \cite{be2} and the references therein). In a direction more
relevant for our development here, using the argument of Ehrhard \cite{er1},
we will show, in section \ref{seceh} below, that (\ref{launo}) follows from
the symmetrization inequality of P\'{o}lya-Szeg\"{o} for Gaussian measure (cf.
\cite{er2} and Section \ref{secprinciple})%
\[
\left\|  \nabla f^{\circ}\right\|  _{L^{2}(\mathbb{R},d\gamma_{1})}%
\leq\left\|  \nabla f\right\|  _{L^{2}(\mathbb{R}^{n},d\gamma_{n})},
\]
where $f^{\circ}$ is the Gaussian symmetric rearrangement of $f$ with respect
to Gaussian measure (cf. Section 2 below).

The purpose of this paper is to give a new approach to lS inequalities through
the use of symmetrization methods. While symmetrization methods are a well
established tool to study Sobolev inequalities, through the combination of
symmetrization and isoperimetric inequalities we uncover new rearrangement
inequalities and connections, that provide a context in which we can treat the
classical and logarithmic Sobolev inequalities in a unified way. Moreover,
with no extra effort we are able to extend the functional lS inequalities to
the general setting of rearrangement invariant spaces. In particular, we
highlight a new extreme embedding which clarifies the connection between lS,
the concentration phenomenon and the John-Nirenberg lemma. Underlying this
last connection is the apparently new observation that concentration
inequalities self improve, a fact we shall treat in detail in a separate paper
(cf. \cite{mami}).

The key to our method are new symmetrization inequalities that involve the
isoperimetric profile and, in this fashion, are strongly associated with
geometric measure theory. In previous papers (cf. \cite{mmp} and the
references therein) we had obtained the corresponding inequalities in the
classical case without making explicit reference to the Euclidean
isoperimetric profile. Using isoperimetry we are able to connect each of the
classical inequalities with their corresponding (new) Gaussian counterparts.
We will show that the difference between the classical and the new Gaussian
inequalities can be simply explained in terms of the difference of the
corresponding isoperimetric profiles. In particular, in the Gaussian case, the
isoperimetric profile is independent of the dimension, and this accounts for
the fact that our rearrangement inequalities in this setting have this
property. Another bonus is that our method is rather general, and amenable to
considerable generalization: to Sobolev inequalities in general measure
spaces, metric Sobolev spaces, even discrete Sobolev spaces. We hope to return
to some of these developments elsewhere.

To describe more precisely our results let us recall that the connection
between isoperimetry and Sobolev inequalities goes back to the work of Maz'ya
and Federer and can be easily explained by combining the formula connecting
the gradient and the perimeter (cf. \cite{ma}):
\begin{equation}
\left\|  \nabla f\right\|  _{1}=\int_{0}^{\infty}Per(\{\left|  f\right|
>t\})dt, \label{burda}%
\end{equation}
with the classical Euclidean isoperimetric inequality:
\begin{equation}
Per(\{\left|  f\right|  >t\})\geq n\varpi_{n}^{1/n}\left(  \left|  \{\left|
f\right|  >t\}\right|  \right)  ^{\frac{n-1}{n}}, \label{burda1}%
\end{equation}
where $\varpi_{n}=$ volume of unit ball in $\mathbb{R}^{n}.$ Indeed, combining
(\ref{burda1}) and (\ref{burda}) yields the sharp form of the
Gagliardo-Nirenberg inequality%
\begin{equation}
(n-1)\varpi_{n}^{1/n}\left\|  f\right\|  _{L^{\frac{n}{n-1},1}(\mathbb{R}%
^{n})}\leq\left\|  \nabla f\right\|  _{L^{1}(\mathbb{R}^{n})}. \label{s5}%
\end{equation}
In \cite{mmp}, we modified Maz'ya's truncation method\footnote{we termed this
method ``symmetrization via truncation''.}, to develop a sharp tool to extract
symmetrization inequalities from Sobolev inequalities like (\ref{s5}). In
particular, we showed that, given any rearrangement invariant norm (r.i. norm)
$\left\|  .\right\|  ,$ the following optimal Sobolev inequality\footnote{This
inequality is optimal and includes the problematic borderline ``end points''
of the $L^{p}$ theory.} holds (cf. \cite{mp})%
\begin{equation}
\left\|  (f^{\ast\ast}(t)-f^{\ast}(t))t^{-1/n}\right\|  \leq c(n,X)\left\|
\nabla f\right\|  ,\text{ }f\in C_{0}^{\infty}(\mathbb{R}^{n}). \label{mp1}%
\end{equation}

An analysis of the role that the power $t^{-1/n}$ plays in this inequality led
us to connect (\ref{mp1}) to isoperimetric profile of $(\mathbb{R}^{n},dx).$
In fact, observe that we can formulate (\ref{burda1}) as%
\[
Per(A)\geq I_{n}(vol_{n}(A)),
\]
where $I_{n}(t)=n\varpi_{n}^{1/n}t^{(n-1)/n}$ is the ``isoperimetric profile''
or the ``isoperimetric function'', and equality is achieved for balls.

The corresponding isoperimetric inequality for Gaussian measure (i.e.
$\mathbb{R}^{n}$ equipped with Gaussian measure $d\gamma_{n}(x)=(2\pi
)^{-n/2}e^{-|x|^{2}/2}dx$), and the solution to the Gaussian isoperimetric
problem, was obtained by Borell \cite{bo} and Sudakov-Tsirelson \cite{sut},
who showed that%
\[
Per(A)\geq I(\gamma_{n}(A)),
\]
with equality archived for half spaces\footnote{In some sense one can consider
half spaces as balls centered at infinity.}, and where $I=$ $I_{\gamma}$ is
the Gaussian profile\footnote{In principle $I$ could depend on $n$ but by the
very definition of half spaces it follows that the Gaussian isoperimetric
profile is dimension free.} (cf. (\ref{abajo}) below for the precise
definition of $I$). To highlight a connection with the lS inequalities, we
only note here that $I$ has the following asymptotic formula near the origin
(say $t\leq1/2$, see Section \ref{Sec2} below),
\begin{equation}
I(t)\simeq t\left(  \log\frac{1}{t}\right)  ^{1/2}. \label{i1}%
\end{equation}

As usual, the symbol $f\simeq g$ will indicate the existence of a universal
constant $c>0$ (independent of all parameters involved) so that $(1/c)f\leq
g\leq c\,f$, while the symbol $f\preceq g$ means that $f\leq c\,g$.

\ 

With this background one may ask: what is the Gaussian replacement of the
Gagliardo-Nirenberg inequality (\ref{s5})? The answer was provided by Ledoux
who showed (cf. \cite{le2})%

\begin{equation}
\int_{0}^{\infty}I(\lambda_{f}(s))ds\leq\int_{\mathbb{R}^{n}}\left|  \nabla
f\right|  d\gamma_{n}(x),\text{ }f\in Lip(\mathbb{R}^{n}). \label{ledox}%
\end{equation}
In fact, following the steps of the proof we indicated for (\ref{s5}), but
using the Gaussian profile instead, we readily arrive at Ledoux's inequality.
This given we were therefore led to apply our method of symmetrization by
truncation to the inequality (\ref{ledox}). We obtained the following
counterpart of (\ref{mp1})
\[
(f^{\ast\ast}(t)-f^{\ast}(t))\leq\frac{t}{I(t)}\left|  \nabla f\right|
^{\ast\ast}(t),
\]
here $f^{\ast}$ denotes the non-increasing rearrangement of $f$ with respec to
to the Lebesgue measure and $f^{\ast\ast}(t)=\frac{1}{t}\int_{0}^{t}f^{\ast
}(s)ds.$ Further analysis showed that, in agreement with the Euclidean case we
had worked out in \cite{mmp}, all these inequalities are in fact
equivalent\footnote{It is somewhat paradoxical that (\ref{launo}), because of
the presence of squares, needs a special treatment and is not, as fas as we
know, equivalent to the isoperimetric inequality.} to the isoperimetric
inequality\footnote{The equivalence between (i) and (ii) in Theorem
\ref{teomain} above is due to Ledoux \cite{le1}.} (cf. Section \ref{sec3} below):

\begin{theorem}
\label{teomain}The following statements are equivalent (all rearrangements are
with respect to Gaussian measure):

(i) Isoperimetric inequality: For every Borel set $A\subset\mathbb{R}^{n},$
with $0<\gamma_{n}(A)<1,$%
\[
Per(A)\geq I(\gamma_{n}(A)).
\]
(ii) Ledoux's inequality: for every Lipschitz function $f$ on $\mathbb{R}%
^{n},$
\begin{equation}
\int_{0}^{\infty}I(\lambda_{f}(s))ds\leq\int_{\mathbb{R}^{n}}\left|  \nabla
f(x)\right|  d\gamma_{n}(x). \label{ledo}%
\end{equation}
(iii) Talenti's inequality (Gaussian version): For every Lipschitz function
$f$ on $\mathbb{R}^{n},$%
\begin{equation}
(-f^{\ast})^{\prime}(s)I(s)\leq\frac{d}{ds}\int_{\{\left|  f\right|  >f^{\ast
}(s)\}}\left|  \nabla f(x)\right|  d\gamma_{n}(x). \label{dosa}%
\end{equation}
(iv) Oscillation inequality (Gaussian version): For every Lipschitz function
$f$ on $\mathbb{R}^{n},$%
\begin{equation}
(f^{\ast\ast}(t)-f^{\ast}(t))\leq\frac{t}{I(t)}\left|  \nabla f\right|
^{\ast\ast}(t). \label{rea}%
\end{equation}
\end{theorem}

This formulation coincides with the corresponding Euclidean result we had
obtained in \cite{mmp}, and thus, in some sense, unifies the classical and
Gaussian Sobolev inequalities. More precisely, by specifying the corresponding
isoperimetric profile we automatically derive the correct results in either
case. Thus, for example, if in (\ref{ledo}) we specify the Euclidean
isoperimetric profile we get the Gagliardo-Nirenberg inequality, in
(\ref{dosa}) we get Talenti's original inequality \cite{Ta} and in (\ref{rea})
we get the rearrangement inequality of \cite{bmr}.

Underlying all these inequalities is the so called P\'{o}lya-Szeg\"{o}
principle. The $L^{p}$ Gaussian versions of this principle had been obtained
earlier by Ehrhard\footnote{For comparison we mention that Ehrhard's results
are formulated in terms of increasing rearrangements.} \cite{er2}. We obtain
here a general version of the P\'{o}lya-Szeg\"{o} principle (cf. \cite{fo}
where the Euclidean case was stated without proof), what may seem surprising
at first is the fact that, in our formulation, the P\'{o}lya-Szeg\"{o}
principle is, in fact, equivalent to the isoperimetric inequality (cf. Section
\ref{secprinciple}).

\begin{theorem}
\label{polya}The following statements are equivalent

(i) Isoperimetric inequality: For every Borel set $A\subset\mathbb{R}^{n},$
with $0<\gamma_{n}(A)<1$%
\[
Per(A)\geq I(\gamma_{n}(A)).
\]
(ii) P\'{o}lya-Szeg\"{o} principle: For every Lipschitz $f$ function on
$\mathbb{R}^{n},$%
\[
\left|  \nabla f^{\circ}\right|  ^{\ast\ast}(s)\leq\left|  \nabla f\right|
^{\ast\ast}(s).
\]
\end{theorem}

Very much like Euclidean symmetrization inequalities lead to optimal Sobolev
and Poincar\'{e} inequalities and embeddings (cf. \cite{mmp}, \cite{mami0} and
the references therein), the new Gaussian counterpart (\ref{rea}) we obtain
here leads to corresponding optimal Gaussian Sobolev-Poincar\'{e} inequalities
as well. The corresponding analog of (\ref{mp1}) is: given any rearrangement
invariant space $X$ on the interval $(0,1),$ we have the optimal inequality,
valid for Lip functions (cf. Section \ref{secpoin} below)%
\begin{equation}
\left\|  f\right\|  _{LS(X)}:=\left\|  \left(  f^{\ast\ast}(t)-f^{\ast
}(t)\right)  \frac{I(t)}{t}\right\|  _{X}\leq\left\|  \nabla f\right\|  _{X}.
\label{corres}%
\end{equation}
The spaces $LS(X)$ defined in this fashion are not necessarily normed,
although often they are equivalent to normed spaces\footnote{For the Euclidean
case a complete study of the normability of these spaces has been recently
given in \cite{pu}.}. As a counterpart to this defect we remark that, since
the Gaussian isoperimetric profile is independent of the dimension, the
inequalities (\ref{corres}) are dimension free. In particular, we note the
following result here (cf. sections \ref{secpoin} and \ref{secfe} below for a
detailed analysis),

\begin{theorem}
\label{opti} Let $X$, $Y$ be two r.i. spaces. Then, the following statements
are equivalent

(i) For every Lipschitz function $f$ on $\mathbb{R}^{n}$
\begin{equation}
\left\|  f-\int f\right\|  _{Y}\preceq\left\|  \nabla f\right\|  _{X}.
\label{poin}%
\end{equation}

(ii) For every positive function $f\in X$ with supp$f\subset(0,1/2),$
\[
\left\|  \int_{t}^{1}f(s)\frac{ds}{I(s)}\right\|  _{Y}\preceq\left\|
f\right\|  _{X}.
\]
Part II. Let $\underline{\alpha}_{X}$ and $\overline{\alpha}_{X}$ be the lower
and the upper Boyd indices of $X$ (see Section \ref{Sec2} below). If
$\underline{\alpha}_{X}>0$, then the following statement is equivalent to (i)
and (ii) above:

(iii)%
\[
\left\|  f\right\|  _{Y}\preceq\left\|  f^{*}(t)\frac{I(t)}{t}\right\|  _{X}.
\]
In particular, if $Y$ is a r.i. space such that (\ref{poin}) holds, then
\[
\left\|  f\right\|  _{Y}\leq\left\|  f^{\ast}(t)\frac{I(t)}{t}\right\|  _{X}.
\]

If $0=\underline{\alpha}_{X}<\overline{\alpha}_{X}<1,$ then the following
statement is equivalent to (i) and (ii) above

(iv)
\[
\left\|  f\right\|  _{Y}\preceq\left\|  f\right\|  _{LS(X)}+\left\|
f\right\|  _{L^{1}.}%
\]
In particular, if $Y$ is a r.i. space such that (\ref{poin}) holds, then
\[
\left\|  f\right\|  _{Y}\leq\left\|  f\right\|  _{LS(X)}+\left\|  f\right\|
_{L^{1}.}.
\]
\end{theorem}

To recognize the logarithmic Sobolev inequalities that are encoded in this
fashion we use the asymptotic property (\ref{i1}) of the isoperimetric profile
$I(s)$ and suitable Hardy type inequalities.

\begin{corollary}
(see Section \ref{secfe} below). Let $X=L^{p},$ $1\leq p<\infty.$ Then,
\[
\int_{0}^{1}\left(  \left(  f-\int f\right)  ^{\ast}(s)\frac{I(s)}{s}\right)
^{p}ds\preceq\int\left|  \nabla f(x)\right|  ^{p}d\gamma_{n}(x).
\]
In particular,%
\[
\int_{0}^{1}f^{\ast}(s)^{p}(\log\frac{1}{s})^{p/2}ds\preceq\int\left|  \nabla
f(x)\right|  ^{p}d\gamma_{n}(x)+\int\left|  f(x)\right|  ^{p}d\gamma_{n}(x).
\]
\end{corollary}

In the final section of this paper we discuss briefly a connection with
concentration inequalities. We refer to Ledoux \cite{le3} for a detailed
account, and detailed references, on the well known connection between lS
inequalities and concentration. In our setting, concentration inequalities can
be derived from a limiting case of the functional lS inequalities. Namely, for
$X=L^{\infty}$ (\ref{corres}) yields%
\[
\left\|  f\right\|  _{LS(L^{\infty})}=\sup_{t<1}\left\{  (f^{\ast\ast
}(t)-f^{\ast}(t))\frac{I(t)}{t}\right\}  \leq\sup_{t}\left|  \nabla f\right|
^{\ast\ast}(t)=\left\|  f\right\|  _{Lip}.
\]
We denote the new space $L_{\log^{1/2}}(\infty,\infty)$ (cf. (\ref{definida})
below). Through the asymptotics of $I(s)$ we see that $L_{\log^{1/2}}%
(\infty,\infty)$ is a variant of the Bennett-DeVore-Sharpley \cite{bds}
space\footnote{$L(\infty,\infty)(R^{n},d\gamma_{n})$ is defined by the
condition%
\[
\sup_{0<t<1}(f^{\ast\ast}(t)-f^{\ast}(t))=\sup_{0<t<1}\frac{1}{t}\int_{0}%
^{t}(f^{\ast}(s)-f^{\ast}(t))ds<\infty.
\]
} $L(\infty,\infty)=$ rearrangement invariant hull of $BMO.$ As it was shown
in \cite{bds}, the definition of $L(\infty,\infty)$ is a reformulation of the
John-Nirenberg inequality and thus yields exponential integrability.
$L_{\log^{1/2}}(\infty,\infty)$ allows us to be more precise about the level
of exponential integrability implied by our inequalities. In this fashion, via
symmetrization and isoperimetry we have connected the John-Nirenberg
inequality with the lS inequalities. 

In a similar manner we can also treat the embedding into $L^{\infty}$ using
the fact that the space $L(\infty,1)=L^{\infty}$ (cf. \cite{bmr}).

Finally, let us state that our main focus in this paper was to develop our
methods and illustrate their reach, but without trying to state the results in
their most general form. We refer the reader to \cite{jmmm} for a general
theory of isoperimetry and symmetrization in the metric setting.

The section headers are self explanatory and provide the organization of the paper.

\section{Gaussian Rearrangements\label{Sec2}}

In this section we review well known results and establish the basic notation
concerning Gaussian rearrangements that we shall use in this paper.

\subsection{Gaussian Profile}

Recall that the $n-$dimensional Gaussian measure on $\mathbb{R}^{n}$ is
defined by
\[
d\gamma_{n}(x)=(2\pi)^{-n/2}e^{-\frac{\left|  x\right|  ^{2}}{2}}%
dx_{1}........dx_{n}.
\]
It is also convenient to let
\[
\phi_{n}(x)=(2\pi)^{-n/2}e^{-\frac{\left|  x\right|  ^{2}}{2}},x\in
\mathbb{R}^{n},
\]
and therefore%
\begin{equation}
\int_{\mathbb{R}^{n}}\phi_{n}(x)dx=\gamma_{n}(\mathbb{R}^{n})=1. \label{proba}%
\end{equation}
Let $\Phi:\mathbb{R}\rightarrow(0,1)$ be the increasing function given by
\[
\Phi(r)=\int_{-\infty}^{r}\phi_{1}(t)dt.
\]
The Gaussian perimeter of a set is defined by
\[
Per(\Omega)=\int_{\partial\Omega}\phi_{n}(x)dH_{n-1}(x),
\]
where $dH_{n-1}(x)$ denotes the Hausdorff $(n-1)$ dimensional measure. The
isoperimetric inequality now reads%
\[
Per(\Omega)\geq I(\gamma_{n}(\Omega)),
\]
where $I$ is the Gaussian isoperimetric function given by (cf. \cite{le1},
\cite{le3})
\begin{equation}
I(t)=\phi_{1}(\Phi^{-1}(t)),\;t\in\lbrack0,1]. \label{abajo}%
\end{equation}

It was shown by Borell \cite{bo} and Sudakov-Tsirelson \cite{sut} that for the
solution of the isoperimetric problem for Gaussian measures we must replace
balls by half spaces. We choose to work with half spaces defined by
\[
H_{r}=\{x=(x_{1},.....x_{n}):x_{1}<r\},\;r\in\mathbb{R}.
\]
Therefore by (\ref{proba}),
\[
\gamma_{n}(H_{r})=\int_{-\infty}^{r}\phi_{1}(t)dt.
\]
Given a measurable set $\Omega\subset\mathbb{R}^{n},$ we let $\Omega^{\circ}$
be the half space defined by
\[
\Omega^{\circ}=H_{r},
\]
where $r\in\mathbb{R}$ is selected so that
\[
\Phi(r)=\gamma_{n}(H_{r})=\gamma_{n}(\Omega).
\]
In other words, $r$ is defined by
\[
r=\Phi^{-1}(\gamma_{n}(\Omega)).
\]
It follows that%
\[
Per(\Omega)\geq Per(\Omega^{\circ})=\phi_{1}(\Phi^{-1}(\gamma_{n}(\Omega))).
\]
Concerning the Gaussian profile $I$ we note here some useful properties for
our development in this paper (cf. \cite{le1} and the references therein).
First, we note that, by direct computation, we have that $I$ satisfies
\begin{equation}
I^{\prime\prime}=\frac{-1}{I}, \label{notada}%
\end{equation}
and, as a consequence of (\ref{proba}), we also have the symmetry
\[
I(t)=I(1-t),t\in\lbrack0,1].
\]
Moreover, from (\ref{notada}) we deduce that $I(s)$ is concave has a maximum
at $t=1/2$ with $I(1/2)=(2\pi)^{-1/2},$ and since $I(0)=0,$ then
$\frac{I(s)-I(0)}{s}$ $=\frac{I(s)}{s}$ is decreasing; summarizing%
\begin{equation}
\frac{I(s)}{s}\text{ is decreasing on }(0,1)\text{ and }\frac{s}{I(s)}\text{
is increasing on }(0,1). \label{ija2}%
\end{equation}

Logarithmic Sobolev inequalities are connected with the asymptotic behavior of
\ $I(t)$ at the origin (or at $1$ by symmetry) (cf. \cite{le1})
\begin{equation}
\lim_{t\rightarrow0}\frac{I(t)}{t(2\log\frac{1}{t})^{1/2}}=1. \label{ija3}%
\end{equation}

\subsection{Rearrangements}

Let $f:\mathbb{R}^{n}\rightarrow\mathbb{R}.$ We define the non increasing,
right continuous, Gaussian distribution function of $f,$ by means of
\[
\lambda_{f}(t)=\gamma_{n}(\{x\in\mathbb{R}^{n}:\left|  f(x)\right|
>t\}),\text{ }t>0.
\]
The rearrangement of $f$ with respect to Gaussian measure, $f^{\ast
}:(0,1]\rightarrow\lbrack0,\infty),$ is then defined, as usual, by
\[
f^{\ast}(s)=\inf\{t\geq0:\lambda_{f}(t)\leq s\},\text{ }t\in(0,1].
\]

In the Gaussian context we replace the classical Euclidean spherical
decreasing rearrangement by a suitable Gaussian substitute, $f^{\circ
}:\mathbb{R}^{n}\rightarrow\mathbb{R}$ defined by
\[
f^{\circ}(x)=f^{\ast}(\Phi(x_{1})).
\]
It is useful to remark here that, as in the Euclidean case, $f^{\circ}$ is
equimeasurable with $f:$%
\begin{align*}
\gamma_{n}(\{x  &  :f^{\circ}(x)>t\})=\gamma_{n}(\{x:f^{\ast}(\Phi
(x_{1}))>t\})\\
&  =\gamma_{n}(\{x:\Phi(x_{1})\leq\lambda_{f}(t)\})\\
&  =\gamma_{n}(\{x:x_{1}\leq\Phi^{-1}(\lambda_{f}(t))\})\\
&  =\gamma_{1}(-\infty,\Phi^{-1}(\lambda_{f}(t)))\\
&  =\lambda_{f}(t).
\end{align*}

\subsection{Rearrangement invariant spaces\label{sec-reor}}

Finally, let us recall briefly the basic definitions and conventions we use
from the theory of rearrangement-invariant (r.i.) spaces, and refer the reader
to \cite{BS} for a complete treatment.

A Banach function space $X=X(\mathbb{R}^{n})$ is called a r.i. space if $g\in
X$ implies that all functions $f$ \ with the same rearrangement with respect
to Gaussian measure, i.e. such that $f^{\ast}=g^{\ast},$ also belong to $X,$
and, moreover, $\Vert f\Vert_{X}=\Vert g\Vert_{X}$. The space $X$ can then be
``reduced'' to a one-dimensional space (which by abuse of notation we still
denote by $X),$ $X=X(0,1),$ consisting of all $g:(0,1)\mapsto\mathbb{R}$ such
that $g^{\ast}(t)=f^{\ast}(t)$ for some function $f\in X$. Typical examples
are the $L^{p}$-spaces and Orlicz spaces.

We shall usually formulate conditions on r.i spaces in terms of the Hardy
operators defined by
\[
Pf(t)=\frac{1}{t}\int_{0}^{t}f(s)ds;\text{ \ \ \ }Qf(t)=\int_{t}^{1}%
f(s)\frac{ds}{s}.
\]
It is well known (see for example \cite[Chapter 3]{BS}), that if $X$ is a r.i.
space, $P$ (resp. $Q$) is bounded on $X$ if and only if the upper Boyd index
$\overline{\alpha}_{X}<1$ (resp. the lower Boyd index $\underline{\alpha}%
_{X}>0$).

We notice for future use that if $X$ is a r.i. space such that $\underline
{\alpha}_{X}>0,$ then the operator
\[
\tilde{Q}f(t)=\left(  1+\log1/t\right)  ^{1/2}\int_{t}^{1}f(s)\frac
{ds}{s\left(  1+\log1/s\right)  ^{1/2}}%
\]
is bounded on $X.$ Indeed, pick $\underline{\alpha}_{X}>a>0,$ then since
$t^{a}\left(  1+\log1/t\right)  ^{1/2}$ is increasing near zero, we get
\[
\tilde{Q}f(t)=\frac{t^{a}\left(  1+\log1/t\right)  ^{1/2}}{t^{a}}\int_{t}%
^{1}f(s)\frac{ds}{s\left(  1+\log1/s\right)  ^{1/2}}\preceq\frac{1}{t^{a}}%
\int_{t}^{1}s^{a}f(s)\frac{ds}{s}=Q_{a}f(t),
\]
and $Q_{a}$ is bounded on $X$ since $\underline{\alpha}_{X}>a$ (see
\cite[Chapter 3]{BS}).

\section{Proof of Theorem \ref{teomain}\label{sec3}}

The proof follows very closely the development in \cite{mmp} with appropriate changes.

$\mathbf{(i)\Rightarrow(ii)}$ By the co-area formula (cf. \cite{ma}) and the
isoperimetric inequality
\begin{align*}
\int\left|  \nabla f(x)\right|  d\gamma_{n}(x)  &  =\int_{0}^{\infty}%
(\int_{\{\left|  f\right|  =s\}}\phi_{n}(x)dH_{n-1}(x))ds\\
&  =\int_{0}^{\infty}Per(\{\left|  f\right|  >s\})ds\\
&  \geq\int_{0}^{\infty}I(\lambda_{f}(s))ds\text{ }.
\end{align*}
$\mathbf{(ii)\Rightarrow(iii)}$ Let $0<t_{1}<t_{2}<\infty.$ The truncations of
$f$ are defined by
\[
f_{t_{1}}^{t_{2}}(x)=\left\{
\begin{array}
[c]{ll}%
t_{2}-t_{1} & \text{if }\left|  f(x)\right|  >t_{2},\\
\left|  f(x)\right|  -t_{1} & \text{if }t_{1}<\left|  f(x)\right|  \leq
t_{2},\\
0 & \text{if }\left|  f(x)\right|  \leq t_{1}.
\end{array}
\right.
\]
Applying (\ref{ledo}) to $f_{t_{1}}^{t_{2}}$ we obtain,
\[
\int_{0}^{\infty}I(\lambda_{f_{t_{1}}^{t_{2}}}(s))ds\leq\int_{\mathbb{R}^{n}%
}\left|  \nabla f_{t_{1}}^{t_{2}}(x)\right|  d\gamma_{n}(x).
\]
We obviously have
\[
\left|  \nabla f_{t_{1}}^{t_{2}}\right|  =\left|  \nabla f\right|
\chi_{\left\{  t_{1}<\left|  f\right|  \leq t_{2}\right\}  },
\]
and, moreover,
\begin{equation}
\int_{0}^{\infty}I(\lambda_{f_{t_{1}}^{t_{2}}}(s))ds=\int_{0}^{t_{2}-t_{1}%
}I(\lambda_{f_{t_{1}}^{t_{2}}}(s))ds. \label{A1}%
\end{equation}
Observe that for $0<s<t_{2}-t_{1}$
\[
\gamma_{n}\left(  \left|  f(x)\right|  \geq t_{2}\right)  \leq\lambda
_{f_{t_{1}}^{t_{2}}}(s)\leq\gamma_{n}\left(  \left|  f(x)\right|
t_{1}\right)  .
\]
Consequently, we have%
\begin{equation}
\int_{0}^{t_{2}-t_{1}}I(\lambda_{f_{t_{1}}^{t_{2}}}(s))ds\geq(t_{2}-t_{1}%
)\min\left(  I(\gamma_{n}\left(  \left|  f\right|  \geq t_{2}\right)
),I(\gamma_{n}\left(  \left|  f\right|  >t_{1}\right)  \right)  . \label{A2}%
\end{equation}
For $s>0$ and $h>0,$ pick $t_{1}=f^{\ast}(s+h),$ $t_{2}=f^{\ast}(s),$ then
\[
s\leq\gamma_{n}\left(  \left|  f(x)\right|  \geq f^{\ast}(s)\right)
\leq\lambda_{f_{t_{1}}^{t_{2}}}(s)\leq\gamma_{n}\left(  \left|  f(x)\right|
f^{\ast}(s+h)\right)  \leq s+h,
\]
Combining (\ref{A1}) and (\ref{A2}) we have,
\begin{align}
\left(  f^{\ast}(s)-f^{\ast}(s+h)\right)  \min(I(s+h),I(s))  &  \leq
\int_{\left\{  f^{\ast}(s+h)<\left|  f\right|  \leq f^{\ast}(s)\right\}
}\left|  \nabla f(x)\right|  d\gamma_{n}(x)\label{truncation}\\
&  \leq\int_{0}^{h}\left|  \nabla f\right|  ^{\ast}(t)dt,\nonumber
\end{align}
whence $f^{\ast}$ is locally absolutely continuous. Thus,
\[
\frac{\left(  f^{\ast}(s)-f^{\ast}(s+h)\right)  }{h}\min(I(s+h),I(s))\leq
\frac{1}{h}\int_{\left\{  f^{\ast}(s+h)<\left|  f\right|  \leq f^{\ast
}(s)\right\}  }\left|  \nabla f(x)\right|  d\gamma_{n}(x).
\]
Letting $h\rightarrow0$ we obtain (\ref{dosa}).\ 

$\mathbf{(iii)\Rightarrow(iv)}$ We will integrate by parts. Let us note first
that using (\ref{truncation}) we have that, for $0<s<t,$
\begin{equation}
s\left(  f^{\ast}(s)-f^{\ast}(t\right)  )\leq\frac{s}{\min(I(s),I(t))}\int
_{0}^{t-s}\left|  \nabla f\right|  ^{\ast}(s)ds. \label{boca}%
\end{equation}
Now,
\begin{align}
f^{\ast\ast}(t)-f^{\ast}(t)  &  =\frac{1}{t}\int_{0}^{t}\left(  f^{\ast
}(s)-f^{\ast}(t)\right)  ds\nonumber\\
&  =\frac{1}{t}\left\{  \left[  s\left(  f^{\ast}(s)-f^{\ast}(t)\right)
\right]  _{0}^{t}+\int_{0}^{t}s\left(  -f^{\ast}\right)  ^{^{\prime}%
}(s)ds\right\} \nonumber\\
&  =\frac{1}{t}\int_{0}^{t}s\left(  -f^{\ast}\right)  ^{^{\prime}%
}(s)ds\nonumber\\
&  =A(t),\nonumber
\end{align}
where the integrated term $\left[  s\left(  f^{\ast}(s)-f^{\ast}(t)\right)
\right]  _{0}^{t}$ vanishes on account of (\ref{boca}). By (\ref{ija2}),
$s/I(s)$ is increasing on $0<s<1$, thus
\begin{align*}
A(t)  &  \leq\frac{1}{I(t)}\int_{0}^{t}I(s)\left(  -f^{\ast}\right)
^{^{\prime}}(s)ds\\
&  \leq\frac{1}{I(t)}\int_{0}^{t}\left(  \frac{\partial}{\partial s}%
\int_{\left\{  \left|  f\right|  >f^{\ast}(s)\right\}  }\left|  \nabla
f(x)\right|  d\gamma_{n}(x)\right)  ds\text{ (by (\ref{dosa}))}\\
&  \leq\frac{1}{I(t)}\int_{\left\{  \left|  f\right|  >f^{\ast}(s)\right\}
}\left|  \nabla f(x)\right|  d\gamma_{n}(x)\text{ }\\
&  \leq\frac{t}{I(t)}\left|  \nabla f\right|  ^{\ast\ast}(t).
\end{align*}
$\mathbf{(iv)\Rightarrow(i)}$ Let $A$ be a Borel set with $0<\gamma_{n}(A)<1.$
We may assume without loss that $Per(A)<\infty.$ By definition we can select a
sequence $\{f_{n}\}_{n\in N}$ of Lip functions such that $f_{n}\underset
{L^{1}}{\rightarrow}\chi_{A}$, and%
\[
Per(A)=\lim\sup_{n\rightarrow\infty}\left\|  \nabla f_{n}\right\|  _{1}.
\]
Therefore,%
\begin{align}
\lim\sup_{n\rightarrow\infty}I(t)(f_{n}^{\ast\ast}(t)-f_{n}^{\ast}(t))  &
\leq\lim\sup_{n\rightarrow\infty}\int_{0}^{t}\left|  \nabla f_{n}(s)\right|
^{\ast}ds\label{bbb}\\
&  \leq\lim\sup_{n\rightarrow\infty}\int\left|  \nabla f_{n}\right|
d\gamma_{n}\nonumber\\
&  =Per(A).\nonumber
\end{align}
As is well known $f_{n}\underset{L^{1}}{\rightarrow}\chi_{A}$ implies that
(cf. \cite[Lemma 2.1]{ga}):
\begin{align*}
f_{n}^{\ast\ast}(t) \rightarrow\chi_{A}^{\ast\ast}(t)  &  , \text{ uniformly
for }t\in\lbrack0,1],\text{and }\\
f_{n}^{\ast}(t)\rightarrow\chi_{A}^{\ast}(t)\text{ }  &  \text{at all points
of continuity of }\chi_{A}^{\ast}.
\end{align*}
Therefore, if we let $r=\gamma_{n}(A),$ and observe that $\chi_{A}^{\ast\ast
}(t)=\min(1,\frac{r}{t}),$ we deduce that for all $t>r,$ $f_{n}^{\ast\ast
}(t)\rightarrow\frac{r}{t},$ and $f_{n}^{\ast}(t)\rightarrow\chi_{A}^{\ast
}(t)=\chi_{(0,r)}(t)=0.$ Inserting this information back in (\ref{bbb}), we
get%
\[
\frac{r}{t}I(t)\leq Per(A),\;\forall t>r.
\]
Now, since $I(t)$ is continuous, we may let $t\rightarrow r$ and we find that
\[
I(\gamma_{n}(A))\leq Per(A),
\]
as we wished to show.

\section{\label{secprinciple}The P\'{o}lya-Szeg\"{o} principle is equivalent
to the isoperimetric inequality}

In this section we prove Theorem \ref{polya}.

Our starting point is inequality (\ref{dosa}). We claim that if $A$ is a
positive Young's function, then
\begin{equation}
A\left(  \left(  -f^{\ast}\right)  ^{^{\prime}}(s)I(s)\right)  \leq
\frac{\partial}{\partial s}\int_{\left\{  \left|  f\right|  >f^{\ast
}(s)\right\}  }A(\left|  \nabla f(x)\right|  )d\gamma_{n}(x). \label{pol}%
\end{equation}
Assuming momentarily the validity of (\ref{pol}), by integration we get%
\begin{equation}
\int_{0}^{1}A\left(  \left(  -f^{\ast}\right)  ^{^{\prime}}(s)I(s)\right)
ds\leq\int_{\mathbb{R}^{n}}A(\left|  \nabla f(x)\right|  )d\gamma_{n}(x).
\label{disfa}%
\end{equation}
It is easy to see that the left hand side is equal to $\int_{\mathbb{R}^{n}%
}A(\left|  \nabla f^{\circ}(x)\right|  )d\gamma_{n}(x).$ Indeed, letting
$s=\Phi(x_{1})$, we find%
\begin{align*}
\int_{0}^{1}A\left(  \left(  -f^{\ast}\right)  ^{^{\prime}}(s)I(s)\right)  ds
&  =\int_{\mathbb{R}}A(\left(  -f^{\ast}\right)  ^{^{\prime}}(\Phi
(x_{1}))I(\Phi(x_{1}))\left|  \Phi^{\prime}(x_{1})\right|  dx\\
&  =\int_{\mathbb{R}^{n}}A(\left(  -f^{\ast}\right)  ^{^{\prime}}(\Phi
(x_{1}))I(\Phi(x_{1}))d\gamma_{n}(x)\\
&  =\int_{\mathbb{R}^{n}}A(\left|  \nabla f^{\circ}(x)\right|  )d\gamma
_{n}(x),
\end{align*}
where in the last step we have used the fact that
\[
\left(  -f^{\ast}\right)  ^{^{\prime}}(\Phi(x_{1}))I(\Phi(x_{1}))=\left(
f^{\ast}\right)  ^{^{\prime}}(\Phi(x_{1}))\Phi^{\prime}(x_{1})=\left|  \nabla
f^{\circ}(x)\right|  .
\]
Consequently, (\ref{disfa}) states that for all Young's functions $A,$ we
have
\[
\int_{\mathbb{R}^{n}}A(\left|  \nabla f^{\circ}(x)\right|  )d\gamma_{n}%
(x)\leq\int_{\mathbb{R}^{n}}A(\left|  \nabla f(x)\right|  )d\gamma_{n}(x),
\]
which, by the Hardy-Littlewood-P\'{o}lya principle, yields
\[
\int_{0}^{t}\left|  \nabla f^{\circ}\right|  ^{\ast}(s)ds\leq\int_{0}%
^{t}\left|  \nabla f\right|  ^{\ast}(s)ds,
\]
as we wished to show.

It remains to prove (\ref{pol}). Here we follow Talenti's argument.\textsl{
}Let $s>0,$ then we have three different alternatives:$\ \mathbf{(a)}$ $s$
belongs to some exceptional set of measure zero, $\mathbf{(b)}$\textsl{\ }%
$\left(  f^{\ast}\right)  ^{^{\prime}}(s)=0,$ or $\mathbf{(c)}$ there is a
neighborhood of $s$ such that $(f^{\ast})^{\prime}(u)$ is not zero, i.e.
$f^{\ast}$ is strictly decreasing. In the two first cases there is nothing to
prove. In case alternative $\mathbf{(c)}$ holds then it follows immediately
from the properties of the rearrangement that for a suitable small $h_{0}>0$
we can write
\[
h=\gamma_{n}\left\{  f^{\ast}(s+h)<\left|  f\right|  \leq f^{\ast}(s)\right\}
,\text{ }0<h<h_{0}.
\]
Therefore, for sufficiently small $h$, we can apply Jensen's inequality to
obtain,
\[
\frac{1}{h}\int_{\left\{  f^{\ast}(s+h)<\left|  f\right|  \leq f^{\ast
}(s)\right\}  }A(\left|  \nabla f(x)\right|  )d\gamma_{n}(x)\geq A\left(
\frac{1}{h}\int_{\left\{  f^{\ast}(s+h)<\left|  f\right|  \leq f^{\ast
}(s)\right\}  }\left|  \nabla f(x)\right|  d\gamma_{n}(x)\right)  .
\]
Arguing like Talenti \cite{Ta} we thus get
\begin{align*}
\frac{\partial}{\partial s}\int_{\left\{  \left|  f\right|  >f^{\ast
}(s)\right\}  }A(\left|  \nabla f(x)\right|  )d\gamma_{n}(x)  &  \geq A\left(
\frac{\partial}{\partial s}\int_{\left\{  \left|  f\right|  >f^{\ast
}(s)\right\}  }\left|  \nabla f(x)\right|  d\gamma_{n}(x)\right) \\
&  \geq A\left(  \left(  -f^{\ast}\right)  ^{^{\prime}}(s)I(s)\right)  ,
\end{align*}
as we wished to show.

To prove the converse we adapt an argument in \cite{bmr}. Let $f$ be a
Lipschitz function on $\mathbb{R}^{n},$ and let $0<t<1$. By the definition of
$f^{\circ}$ we can write%
\begin{align*}
f^{\ast}(t)-f^{\ast}(1^{-})  &  =f^{\ast}(\Phi(\Phi^{-1}(t)))-f^{\ast}%
(\Phi(\infty))\\
&  =\int_{\Phi^{-1}(t)}^{\infty}\left|  \nabla f^{\circ}\right|  (s)ds.
\end{align*}
Thus,%
\[
f^{\ast\ast}(t)-f^{\ast}(1^{-})=\frac{1}{t}\int_{0}^{t}\int_{\Phi^{-1}%
(r)}^{\infty}\left|  \nabla f^{\circ}\right|  (s)dsdr.
\]
Making the change of variables $s=\Phi^{-1}(z)$ in the inner integral and then
changing the order of integration$,$ we find
\begin{align*}
f^{\ast\ast}(t)-f^{\ast}(1^{-})  &  =\frac{1}{t}\int_{0}^{t}\int_{r}%
^{1}\left|  \nabla f^{\circ}\right|  (\Phi^{-1}(z))\left(  \Phi^{-1}%
(z)\right)  ^{^{\prime}}dzdr\\
&  =\int_{t}^{1}\left|  \nabla f^{\circ}\right|  (\Phi^{-1}(z))\left(
\Phi^{-1}(z)\right)  ^{^{\prime}}dz+\frac{1}{t}\int_{0}^{t}z\left|  \nabla
f^{\circ}\right|  (\Phi^{-1}(z))\left(  \Phi^{-1}(z)\right)  ^{^{\prime}}dz\\
&  =f^{\ast}(t)-f^{\ast}(1^{-})+\frac{1}{t}\int_{0}^{t}z\left|  \nabla
f^{\circ}\right|  (\Phi^{-1}(z))\left(  \Phi^{-1}(z)\right)  ^{^{\prime}}dz.
\end{align*}
Since $\Phi^{\prime}(\Phi^{-1}(z))=\phi_{1}(\Phi^{-1}(z)))=I(z),$ we readily
deduce that $\left(  \Phi^{-1}(z)\right)  ^{^{\prime}}=\frac{1}{I(z)}.$ Thus,
\[
f^{\ast\ast}(t)-f^{\ast}(1^{-})=f^{\ast}(t)-f^{\ast}(1^{-})+\frac{1}{t}%
\int_{0}^{t}z\left|  \nabla f^{\circ}\right|  (\Phi^{-1}(z))\frac{1}{I(z)}dz,
\]
and consequently%
\begin{align*}
f^{\ast\ast}(t)-f^{\ast}(t)  &  =\frac{1}{t}\int_{0}^{t}z\left|  \nabla
f^{\circ}\right|  (\Phi^{-1}(z))\frac{1}{I(z)}dz\\
&  \leq\frac{t}{I(t)}\frac{1}{t}\int_{0}^{t}z\left|  \nabla f^{\circ}\right|
(\Phi^{-1}(z))dz\text{ (since }t/I(t)\text{ is increasing)}\\
&  =\frac{1}{I(t)}\int_{-\infty}^{\Phi^{-1}(t)}\left|  \nabla f^{\circ
}\right|  (s)\Phi^{\prime}(s)ds\\
&  =\frac{1}{I(t)}\int_{-\infty}^{\Phi^{-1}(t)}\left|  \nabla f^{\circ
}\right|  (s)d\gamma_{1}(s)\\
&  \leq\int_{0}^{t}\left|  \nabla f^{\circ}\right|  ^{\ast}(s)ds\text{ (since
}\gamma_{1}(-\infty,\Phi^{-1}(t))=t\text{).}%
\end{align*}
Summarizing, we have shown that
\[
(f^{\ast\ast}(t)-f^{\ast}(t))\leq\frac{t}{I(t)}\left|  \nabla f^{\circ
}\right|  ^{\ast\ast}(t),
\]
which combined with our current hypothesis yields
\[
(f^{\ast\ast}(t)-f^{\ast}(t))\leq\frac{t}{I(t)}\left|  \nabla f^{\circ
}\right|  ^{\ast\ast}(t)\leq\frac{t}{I(t)}\left|  \nabla f\right|  ^{\ast\ast
}(t).
\]
By Theorem \ref{teomain} the last inequality is equivalent to the
isoperimetric inequality.

\begin{remark}
We note here, for future use, that the discussion in this section shows that
the following equivalent form of the P\'{o}lya-Szeg\"{o} principle holds
\[
\int_{0}^{t}((-f^{\ast})^{\prime}(.)I(.))^{\ast}(s)ds\leq\int_{0}^{t}\left|
\nabla f\right|  ^{\ast}(s)ds.
\]
Therefore, by the Hardy-Littlewood principle, for every r.i. space $X$ on
$(0,1),$%
\begin{equation}
\left\|  \left(  -f^{\ast}\right)  ^{^{\prime}}(s)I(s)\right\|  _{X}%
\leq\left\|  \nabla f\right\|  _{X}. \label{provadas}%
\end{equation}
\end{remark}

\section{\label{seceh}The P\'{o}lya-Szeg\"{o} principle implies Gross' inequality}

We present a proof due to Ehrhard \cite{er1}, showing that the P\'{o}%
lya-Szeg\"{o} principle implies (\ref{launo}). We present full details, since
Ehrhard's method is apparently not well known and some details are missing in
\cite{er1}.

We first prove a one dimensional inequality which, by symmetrization and
tensorization, will lead to the desired result.

Let $f:\mathbb{R}\rightarrow\mathbb{R}$ be a Lip function such that $f$ and
$f^{\prime}\in L^{1}.$ By Jensen's inequality%
\begin{align*}
\int_{-\infty}^{\infty}\left|  f(x)\right|  \ln\left|  f(x)\right|  dx  &
=\left\|  f\right\|  _{L^{1}}\int_{-\infty}^{\infty}\ln\left|  f(x)\right|
\frac{\left|  f(x)\right|  dx}{\left\|  f\right\|  _{L^{1}}}\\
&  \leq\left\|  f\right\|  _{L^{1}}\ln(\int_{-\infty}^{\infty}\left|
f(x)\right|  \frac{\left|  f(x)\right|  dx}{\left\|  f\right\|  _{L^{1}}}).
\end{align*}
We estimate the inner integral using the fundamental theorem of Calculus:
$\left|  f(x)\right|  \leq\left\|  f^{\prime}\right\|  _{L^{1}},$ to obtain%
\[
\int_{-\infty}^{\infty}\left|  f(x)\right|  \ln\left|  f(x)\right|
dx\leq\left\|  f\right\|  _{L^{1}}\ln\left\|  f^{\prime}\right\|  _{L^{1}}.
\]
Applying the preceding to $f^{2}$ we get:%
\[
\int_{-\infty}^{\infty}\left|  f(x)\right|  ^{2}\ln\left|  f(x)\right|
dx\leq\frac{1}{2}\left\|  f\right\|  _{L^{2}}^{2}\ln2\left\|  ff^{^{\prime}%
}\right\|  _{L^{1}}.
\]
Using H\"{o}lder's inequality $\left\|  ff^{^{\prime}}\right\|  _{L^{1}}%
\leq\left\|  f\right\|  _{L^{2}}\left\|  f^{^{\prime}}\right\|  _{L^{2}},$ and
elementary properties of the logarithm we find%
\begin{align}
\int_{-\infty}^{\infty}\left|  f(x)\right|  ^{2}\ln\left|  f(x)\right|  dx  &
\leq\frac{1}{2}\left\|  f\right\|  _{L^{2}}^{2}\ln2\left\|  f\right\|
_{L^{2}}\left\|  f^{^{\prime}}\right\|  _{L^{2}}\label{level4}\\
&  =\frac{1}{4}\left\|  f\right\|  _{L^{2}}^{2}\ln4\left\|  f\right\|
_{L^{2}}^{4}\frac{\left\|  f^{^{\prime}}\right\|  _{L^{2}}^{2}}{\left\|
f\right\|  _{L^{2}}^{2}}\nonumber\\
&  =\frac{1}{4}\left\|  f\right\|  _{L^{2}}^{2}\ln4\frac{\left\|  f^{^{\prime
}}\right\|  _{L^{2}}^{2}}{\left\|  f\right\|  _{L^{2}}^{2}}+\left\|
f\right\|  _{L^{2}}^{2}\ln\left\|  f\right\|  _{L^{2}}\nonumber\\
&  \leq\left\|  f^{^{\prime}}\right\|  _{L^{2}}^{2}+\left\|  f\right\|
_{L^{2}}^{2}\ln\left\|  f\right\|  _{L^{2}}\text{ (in the last step we used
}\ln t\leq t).\nonumber
\end{align}
We apply (\ref{level4}) to $u=(2\pi e^{x^{2}})^{-1/4}f(x)=\phi_{1}%
(x)^{1/2}f(x)$ and compute both sides of (\ref{level4})$.$ The left hand side
becomes%
\begin{align*}
\int_{-\infty}^{\infty}\left|  u(x)\right|  ^{2}\ln\left|  u(x)\right|  dx  &
=\int_{-\infty}^{\infty}\left|  f(x)\right|  ^{2}\left(  \ln\left|
f(x)\right|  +\ln(2\pi e^{x^{2}})^{-1/4}\right)  d\gamma_{1}(x)\\
&  =\int_{-\infty}^{\infty}\left|  f(x)\right|  ^{2}\ln\left|  f(x)\right|
d\gamma_{1}(x)-\frac{1}{4}\ln2\pi\left\|  f\right\|  _{L^{2}(d\gamma_{1})}\\
&  -\frac{1}{4}\int_{-\infty}^{\infty}\left|  f(x)\right|  ^{2}x^{2}%
d\gamma_{1}(x),
\end{align*}
while the right hand side is equal to%
\begin{align}
\left\|  f^{^{\prime}}\right\|  _{L^{2}}^{2}  &  =\left\|  f^{^{\prime}%
}\right\|  _{L^{2}(d\gamma_{1})}^{2}+\frac{1}{4}\int_{-\infty}^{\infty
}f(x)^{2}x^{2}d\gamma_{1}(x)-\int_{-\infty}^{\infty}f^{^{\prime}}%
(x)f(x)x\phi_{1}(x)dx\label{level2}\\
&  =\left\|  f^{^{\prime}}\right\|  _{L^{2}(d\gamma_{1})}^{2}-\frac{1}{4}%
\int_{-\infty}^{\infty}f(x)^{2}x^{2}d\gamma_{1}(x)+\frac{1}{2}\int_{-\infty
}^{\infty}f(x)^{2}x^{2}d\gamma_{1}(x)\nonumber\\
&  -\int_{-\infty}^{\infty}f^{^{\prime}}(x)f(x)x\phi_{1}(x)dx.\nonumber
\end{align}
We simplify the last expression integrating by parts the third integral to the
right,%
\begin{align*}
\frac{1}{2}\int_{-\infty}^{\infty}f(x)^{2}x^{2}d\gamma_{1}(x)  &  =-\frac
{1}{2}\int_{-\infty}^{\infty}f(x)^{2}xd(((2\pi)^{-1/2}e^{-x^{2}}))\\
&  =\left.  -\frac{1}{2}f(x)^{2}x((2\pi)^{-1/2}e^{-x^{2}})\right|  _{-\infty
}^{\infty}+\\
&  \frac{1}{2}\int_{-\infty}^{\infty}((2\pi)^{-1/2}e^{-x^{2}})[2f(x)f^{\prime
}(x)x+f^{2}(x)]dx\\
&  =\int_{-\infty}^{\infty}f(x)f^{\prime}(x)x\phi_{1}(x)dx+\frac{1}{2}\left\|
f\right\|  _{L^{2}(d\gamma_{1})}^{2}.
\end{align*}
We insert this back in (\ref{level2}) and then comparing results and
simplifying we arrive at%
\begin{align}
\int_{-\infty}^{\infty}\left|  f(x)\right|  ^{2}\ln\left|  f(x)\right|
d\gamma_{1}(x)  &  \leq\left\|  f^{^{\prime}}\right\|  _{L^{2}(d\gamma_{1}%
)}^{2}+\left\|  f\right\|  _{L^{2}(d\gamma_{1})}^{2}\ln\left\|  f\right\|
_{L^{2}(d\gamma_{1})}^{2}\label{estuario}\\
&  + \frac{\ln(2\pi e^{2})}{4}\left\|  f\right\|  _{L^{2}(d\gamma_{1})}%
^{2}.\nonumber
\end{align}

Let $f$ be a Lipchitz function on $\mathbb{R}^{n}$. We form the symmetric
rearrangement $f^{\circ}$ considered as a one dimensional function. Then,
(\ref{estuario}) applied to $f^{\circ}$, combined with the fact that
$f^{\circ}$ is equimesurable with $f$ and the P\'{o}lya-Szeg\"{o} principle,
yields%
\begin{align}
\int_{\mathbb{R}^{n}}\left|  f(x)\right|  ^{2}\ln\left|  f(x)\right|
d\gamma_{n}(x)  &  =\int_{\mathbb{R}}\left|  f^{\circ}(x)\right|  ^{2}%
\ln\left|  f^{\circ}(x)\right|  d\gamma_{1}(x)\label{levelr}\\
&  \leq\left\|  f^{^{\circ\prime}}\right\|  _{L^{2}(d\gamma_{1})}^{2}+\left\|
f\right\|  _{L^{2}(d\gamma_{1})}^{2}\ln\left\|  f\right\|  _{L^{2}(d\gamma
_{1})}^{2}\nonumber\\
&  + \frac{\ln(2\pi e^{2})}{4}\left\|  f\right\|  _{L^{2}(d\gamma_{1})}%
^{2}\nonumber\\
&  =\left\|  \left|  \nabla f^{\circ}(x)\right|  \right\|  _{L^{2}(d\gamma
_{n})}^{2}+\left\|  f^{\circ}\right\|  _{L^{2}(d\gamma_{n})}^{2}\ln\left\|
f^{\circ}\right\|  _{L^{2}(d\gamma_{n})}^{2}\nonumber\\
&  + \frac{\ln(2\pi e^{2})}{4}\left\|  f^{\circ}\right\|  _{L^{2}(d\gamma
_{n})}^{2}\nonumber\\
&  \leq\left\|  \nabla f\right\|  _{L^{2}(d\gamma_{n})}^{2}+\left\|
f\right\|  _{L^{2}(d\gamma_{n})}^{2}\ln\left\|  f\right\|  _{L^{2}(d\gamma
_{n})}^{2}\nonumber\\
&  + \frac{\ln(2\pi e^{2})}{4}\left\|  f\right\|  _{L^{2}(d\gamma_{n})}%
^{2}.\nonumber
\end{align}

We now use tensorization to prove (\ref{launo}). Note that, by homogeneity, we
may assume that $f$ has been normalized so that $\left\|  f\right\|
_{L^{2}(d\gamma_{n})}=1.$ Let $l\in N,$ and let $F$ be defined on
$(\mathbb{R}^{n})^{l}=\mathbb{R}^{nl}$ by $F(x)=\prod_{k=1}^{l}f(x_{k}),$
where $x_{k}\in\mathbb{R}^{n},k=1,..l.$ The $\mathbb{R}^{nl}$ version of
(\ref{levelr}) applied to $F,$ and translated back in terms of $f,$ yields%
\[
l\int_{\mathbb{R}^{n}}\left|  f(x)\right|  ^{2}\ln\left|  f(x)\right|
d\gamma_{n}(x)\leq l\left\|  \nabla f\right\|  _{L^{2}(d\gamma_{n})}^{2}%
+\frac{\ln(2\pi e^{2})}{4}.
\]
Therefore, upon diving by $l$ and letting $l\rightarrow\infty,$ we obtain%
\[
\int_{\mathbb{R}^{n}}\left|  f(x)\right|  ^{2}\ln\left|  f(x)\right|
d\gamma_{n}(x)\leq\left\|  \nabla f\right\|  _{L^{2}(d\gamma_{n})}^{2},
\]
as we wished to show.

\section{Poincar\'{e} type inequalities\label{secpoin}}

We consider $L^{1}$ Poincar\'{e} inequalities first. Indeed, for $L^{1}$ norms
the Poincar\'{e} inequalities are a simple variant of Ledoux's inequality. Let
$f$ be a Lipschitz function on $\mathbb{R}^{n},$ and let $m$ a
median\footnote{i.e. $\gamma_{n}\left(  f\geq m\right)  \geq1/2$ and
$\gamma_{n}\left(  f\leq m\right)  \geq1/2$.} of $f.$ Set $f^{+}=\max(f-m,0)$
and $f^{-}=-\min(f-m,0)$ so that $f-m=f^{+}-f^{-}.$ Then,%
\begin{align*}
\int_{\mathbb{R}^{n}}\left|  f-m\right|  d\gamma_{n}  &  =\int_{\mathbb{R}%
^{n}}f^{+}d\gamma_{n}+\int_{\mathbb{R}^{n}}f^{-}d\gamma_{n}\\
&  =\int_{0}^{\infty}\lambda_{f^{+}}(s)ds+\int_{0}^{\infty}\lambda_{f^{-}%
}(s)ds\\
&  =(A)
\end{align*}
We estimate each of these integrals using the properties of the isoperimetric
profile and Ledoux's inequality (\ref{ledo}). First we use the fact that
$\frac{I(s)}{s}$ is decreasing on $0<s<1/2,$ combined with the definition of
median, to find that%
\[
2\lambda_{g}(s)I(\frac{1}{2})\leq I(\lambda_{g}(s)),\text{ where }%
g=f^{+}\text{ or }g=f^{-}.
\]
Consequently,%
\begin{align*}
(A)  &  \leq\frac{1}{2I(\frac{1}{2})}\left(  \int_{0}^{\infty}I(\lambda
_{f^{+}}(s))ds+\int_{0}^{\infty}I(\lambda_{f^{-}}(s))ds\right) \\
&  \leq\frac{1}{2I(\frac{1}{2})}\left(  \int_{\mathbb{R}^{n}}\nabla
f^{+}(x)d\gamma_{n}(x)+\int_{\mathbb{R}^{n}}\nabla f^{+}(x)d\gamma
_{n}(x)\right)  \text{ (by (\ref{ledo})}\\
&  =\frac{1}{2I(1/2)}\int_{\mathbb{R}^{n}}\left|  \nabla f(x)\right|
d\gamma_{n}(x).
\end{align*}
Thus,%
\begin{equation}
\int_{\mathbb{R}^{n}}\left|  f-m\right|  d\gamma_{n}\leq\frac{1}{2I(1/2)}%
\int_{\mathbb{R}^{n}}\left|  \nabla f(x)\right|  d\gamma_{n}(x).
\label{poinca1}%
\end{equation}

\textbf{We now prove Theorem \ref{opti}}.

\begin{proof}
$\mathbf{(i)\rightarrow(ii).}$ Obviously condition (\ref{poin}) is equivalent
to
\[
\left\|  f-m\right\|  _{Y}\preceq\left\|  \nabla f\right\|  _{X},
\]
where $m$ is a median of $f.$ Let $f$ be a\ positive measurable function with
$suppf$ $\subset(0,1/2).$ Define%
\[
u(x)=\int_{\Phi(x_{1})}^{1}f(s)\frac{ds}{I(s)},\text{ \ \ }x\in\mathbb{R}%
^{n}.
\]
It is plain that $u$ is a Lipschitz function on $\mathbb{R}^{n}$ such
that$\ \gamma_{n}\left(  u=0\right)  \geq1/2,$ and therefore it has $0$
median. Moreover,
\[
\left|  \nabla u(x)\right|  =\left|  \frac{\partial}{\partial x_{1}%
}u(x)\right|  =\left|  -f(\Phi(x_{1}))\frac{\Phi^{\prime}(x_{1})}{I(\Phi
(x_{1}))}\right|  =f(\Phi(x_{1})).
\]
It follows that%
\[
u^{\ast}(t)=\int_{t}^{1}f(s)\frac{ds}{I(s)},\text{ and }\left|  \nabla
u\right|  ^{\ast}(t)=f^{\ast}(t).
\]
Consequently, from%
\[
\left\|  u-0\right\|  _{Y}\preceq\left\|  \nabla u\right\|  _{X}%
\]
we deduce that%
\[
\left\|  \int_{t}^{1}f(s)\frac{ds}{I(s)}\right\|  _{Y}\preceq\left\|
f\right\|  _{X}.
\]
$\mathbf{(ii)\rightarrow(i).}$ Let $f$ be a Lipschitz function $f$ on
$\mathbb{R}^{n}$. Write%
\[
f^{\ast}(t)=\int_{t}^{1/2}\left(  -f^{\ast}\right)  ^{^{\prime}}(s)ds+f^{\ast
}(1/2).
\]
Thus,
\begin{align*}
\left\|  f\right\|  _{Y}  &  =\left\|  f^{\ast}\right\|  _{Y}\leq2\left\|
f^{\ast}\chi_{\lbrack0,1/2]}\right\|  _{Y}\preceq\left\|  \int_{t}%
^{1/2}\left(  -f^{\ast}\right)  ^{^{\prime}}(s)ds\right\|  _{Y}+f^{\ast
}(1/2)\left\|  1\right\|  _{Y}\\
&  \leq\left\|  \int_{t}^{1/2}\left(  -f^{\ast}\right)  ^{^{\prime}%
}(s)I(s)\frac{ds}{I(s)}\right\|  _{Y}+2\left\|  1\right\|  _{Y}\left\|
f\right\|  _{L_{1}}\\
&  \preceq\left\|  \left(  -f^{\ast}\right)  ^{^{\prime}}(s)I(s)\right\|
_{X}+\left\|  f\right\|  _{L_{1}}\\
&  \preceq\left\|  \nabla f\right\|  _{X}\text{ (by (\ref{poinca1}) and
(\ref{provadas})).}%
\end{align*}

\textbf{Part II}. Case $0<$\underline{$\alpha$}$_{X}:$

$\mathbf{(ii)\rightarrow(iii)}$ Let $0<t<1/4$, then%
\[
f^{\ast}(2t)\preceq\int_{t}^{2t}f^{\ast}(s)\frac{ds}{s}\leq\int_{t}%
^{1/2}f^{\ast}(s)\frac{I(s)}{s}\frac{ds}{I(s)},
\]
therefore,%
\begin{align*}
\left\|  f^{\ast}(2t)\right\|  _{Y}  &  \preceq\left\|  \int_{t}^{1/2}f^{\ast
}(s)\frac{I(s)}{s}\frac{ds}{I(s)}\right\|  _{Y}+f^{\ast}(1/2)\\
&  \preceq\left\|  f^{\ast}(t)\frac{I(t)}{t}\right\|  _{X}+f^{\ast}(1/2)\text{
\ (by (ii))}\\
&  \preceq\left\|  f^{\ast}(t)\frac{I(t)}{t}\right\|  _{X}+\left\|  f\right\|
_{1}\\
&  \preceq\left\|  f^{\ast}(t)\frac{I(t)}{t}\right\|  _{X}.
\end{align*}
$\mathbf{(iii)\rightarrow(ii)}$ By hypothesis%
\[
\left\|  \int_{t}^{1/2}f^{\ast}(s)\frac{ds}{I(s)}\right\|  _{Y}\preceq\left\|
\left(  \int_{t}^{1/2}f^{\ast}(s)\frac{ds}{I(s)}\right)  \frac{I(t)}%
{t}\right\|  _{X}.
\]
Using that (see \ref{ija3}),
\[
\frac{I(s)}{s}\simeq\sqrt{\log\frac{1}{s}}\simeq\sqrt{1+\log\frac{1}{s}%
},\text{ \ \ \ }0<s<1/2
\]
we have
\[
\left(  \int_{t}^{1/2}f(s)\frac{ds}{I(s)}\right)  \frac{I(t)}{t}\preceq
\sqrt{1+\log\frac{1}{t}}\int_{t}^{1}f(s)\frac{ds}{s\sqrt{1+\log\frac{1}{s}}%
}=\tilde{Q}f(t).
\]
Now. from \underline{$\alpha$}$_{X}>0$ it follows that $\tilde{Q}$ is a
bounded operator on $X$ (see Section \ref{sec-reor}) and thus we are able to conclude.

\textbf{Part II}. Case $0=$\underline{$\alpha$}$_{X}<\overline{\alpha}_{X}<1:$

$\mathbf{(ii)\rightarrow(iv)}$ By the fundamental theorem of Calculus and
(ii), we have%
\begin{align*}
\left\|  f^{\ast\ast}\chi_{(0,1/2)}\right\|  _{Y}  &  \preceq\left\|  \int
_{t}^{1/2}\left(  f^{\ast\ast}(s)-f^{\ast}(s)\right)  \frac{ds}{s}\right\|
_{Y}+f^{\ast\ast}(1/2)\left\|  1\right\|  _{Y}\\
&  \preceq\left\|  \int_{t}^{1}\frac{I(s)}{s}\left(  f^{\ast\ast}(s)-f^{\ast
}(s)\right)  \chi_{(0,1/2)}(s)\frac{ds}{I(s)}\right\|  _{Y}+\left\|
f\right\|  _{1}\\
&  \preceq\left\|  (f^{\ast\ast}(t)-f^{\ast}(t))\chi_{(0,1/2)}(t)\frac
{I(t)}{t}\right\|  _{X}+\left\|  f\right\|  _{1}\\
&  \preceq\left\|  (f^{\ast\ast}(t)-f^{\ast}(t))\frac{I(t)}{t}\right\|
_{X}+\left\|  f\right\|  _{1}.
\end{align*}

$\mathbf{(iv)\rightarrow(i)}$ Let $f$ be a Lipschitz function on
$\mathbb{R}^{n}$, let $m$ be a median of $f$ and let $g=f-m.$ By hypothesis we
have
\[
\left\|  g\right\|  _{Y}\preceq\left\|  (g^{\ast\ast}(t)-g^{\ast}%
(t))\frac{I(t)}{t}\right\|  _{X}+\left\|  g\right\|  _{1}.
\]

From\ (see \cite{bmr})%
\[
g^{\ast\ast}(t)-g^{\ast}(t)\leq P(g^{\ast}(s/2)-g^{\ast}(s))(t)+g^{\ast
}(t/2)-g^{\ast}(t),
\]
and using the fact that $\frac{I(t)}{t}$ decreases,
\[
P(g^{\ast}(s/2)-g^{\ast}(s))(t)\frac{I(t)}{t}\leq P(g^{\ast}(s/2)-g^{\ast
}(s)\frac{I(s)}{s})(t).
\]
Therefore,
\begin{align*}
\left\|  (g^{\ast\ast}(t)-g^{\ast}(t))\frac{I(t)}{t}\right\|  _{X}  &
\leq\left\|  P(g^{\ast}(s/2)-g^{\ast}(s)\frac{I(s)}{s})(t)\right\|
_{X}+\left\|  (g^{\ast}(t/2)-g^{\ast}(t)\frac{I(t)}{t}\right\|  _{X}\\
&  \preceq\left\|  (g^{\ast}(t/2)-g^{\ast}(t)\frac{I(t)}{t}\right\|
_{X}\text{ \ \ (since }\overline{\alpha}_{X}<1\text{).}%
\end{align*}
We compute the right hand side,%
\begin{align*}
\left\|  (g^{\ast}(t/2)-g^{\ast}(t)\frac{I(t)}{t}\right\|  _{X}  &  =\left\|
\left(  \int_{t/2}^{t}\left(  -g^{\ast}\right)  ^{^{\prime}}(s)ds\right)
\frac{I(t)}{t}\right\|  _{X}\\
&  \leq\left\|  \int_{t/2}^{t}\left(  -g^{\ast}\right)  ^{^{\prime}}%
(s)\frac{I(s)}{s}ds\right\|  _{X}\\
&  \leq\left\|  \frac{2}{t}\int_{t/2}^{t}\left(  -g^{\ast}\right)  ^{^{\prime
}}(s)I(s)ds\right\|  _{X}\\
&  \leq2\left\|  \frac{1}{t}\int_{0}^{t}\left(  -g^{\ast}\right)  ^{^{\prime}%
}(s)I(s)ds\right\|  _{X}\\
&  \preceq\left\|  \left(  -g^{\ast}\right)  ^{^{\prime}}(t)I(t)\right\|
_{X}\\
&  \preceq\left\|  \nabla f\right\|  _{X}\text{ \ \ (by (\ref{provadas})).}%
\end{align*}
Summarizing, we have obtained%
\[
\left\|  g\right\|  _{Y}\preceq\left\|  \nabla f\right\|  _{X}+\left\|
g\right\|  _{1}\preceq\left\|  \nabla f\right\|  _{X}\text{ \ \ (by
(\ref{poinca1})).}%
\]
\end{proof}

\subsection{Feissner type inequalities\label{secfe}}

Theorem \ref{opti}) readily implies Feissner's inequalities (\ref{launodos}).
Indeed, for the particular choice $X=L^{p}$ ($1\leq p<\infty)$, Theorem
\ref{opti} yields
\[
\int_{0}^{1}\left(  \left(  f-\int f\right)  ^{\ast}(s)\frac{I(s)}{s}\right)
^{p}ds\preceq\int\left|  \nabla f(x)\right|  ^{p}d\gamma_{n}(x).
\]
In particular, using again the asymptics of $I(s),$ $0<s<1/2$, we get
\[
\int_{0}^{1}f^{\ast}(s)^{p}(\log\frac{1}{s})^{p/2}ds\preceq\int\left|  \nabla
f(x)\right|  ^{p}d\gamma_{n}(x)+\int\left|  f(x)\right|  ^{p}d\gamma_{n}(x).
\]
Moreover, the space $L^{p}(LogL)^{1/2}$ is best possible among r.i. spaces $Y$
for which the Poincar\'{e} inequality $\Vert f-\int f\Vert_{Y}\preceq
\Vert\nabla f\Vert_{L^{p}}$ holds.

The case $X=L^{\infty},$ which is new is more interesting. Indeed, since
$I(t)/t$ decreases,
\[
\sup_{0<t<1}f^{\ast}(t)\frac{I(t)}{t}<\infty\text{ }\Leftrightarrow\text{
}f=0.
\]
But Theorem \ref{opti} ensures that
\begin{equation}
\left\|  \left(  \left(  f-\int f\right)  ^{\ast\ast}(t)-\left(  f-\int
f\right)  ^{\ast}(t)\right)  \frac{I(t)}{t}\right\|  _{L^{\infty}}%
\preceq\left\|  \nabla f\right\|  _{L^{\infty}}. \label{embl1}%
\end{equation}
Furthermore, for every r.i space $Y$ such that%
\[
\left\|  f-\int f\right\|  _{Y}\preceq\left\|  \nabla f\right\|  _{L^{\infty}%
},
\]
the following embedding holds
\[
\left\|  f\right\|  _{Y}\preceq\left\|  \left(  f^{\ast\ast}(t)-f^{\ast
}(t)\right)  \frac{I(t)}{t}\right\|  _{L^{\infty}}+\left\|  f\right\|  _{1}.
\]
Notice that due to the cancellation afforded by $f^{\ast\ast}(t)-f^{\ast}(t)$,
the corresponding space $LS(L^{\infty})$ is nontrivial. The relation between
concentration and $LS(L^{\infty})$ will be studied in the next section.

\section{On limiting embeddings and concentration\label{secconc}}

Elsewhere \footnote{In particular the method of symmetrization by truncation
can be extended to this setting.} (cf. \cite{mami}) we shall explore in detail
the connection between concentration inequalities and symmetrization,
including the self improving properties of concentration. In this section we
merely wish to call attention to the connection between a limiting lS
inequality that follows from (\ref{rea}) and concentration. We have argued
that, in \ the Gaussian world, Ledoux's embedding corresponds to the
Gagliardo-Nirenberg embedding. In the classical $n-$dimensional Euclidean case
the ``other'' borderline case for the Sobolev embedding theorem occurs when
the index of integrability of the gradients in the Sobolev space, say $p,$ is
equal to the dimension i.e. $p=n.$ In this case, as is well known, from
$\left|  \nabla f\right|  \in L^{n}(\mathbb{R}^{n})$ we can deduce the
exponential integrability of $\left|  f\right|  ^{n^{\prime}}$ (cf.
\cite{tr})$.$ A refinement of this result, which follows from the Euclidean
version of (\ref{rea}), is given by the following inequality from \cite{bmr}%
\[
\left\{  \int_{0}^{\infty}(f^{\ast\ast}(s)-f^{\ast}(s))^{n}\frac{ds}%
{s}\right\}  ^{1/n}\preceq\left\{  \int_{0}^{\infty}\left|  \nabla
f(x)\right|  ^{n}dx\right\}  ^{1/n}.
\]
In this fashion one could consider the corresponding borderline Gaussian
embedding that results from (\ref{rea}) when $n=p=\infty.$ The result now
reads%
\begin{equation}
\sup_{t<1}\left\{  (f^{\ast\ast}(t)-f^{\ast}(t))\frac{I(t)}{t}\right\}
\leq\sup_{t}\left|  \nabla f\right|  ^{\ast\ast}(t)=\left\|  f\right\|
_{Lip}. \label{caso-inf}%
\end{equation}
We now show how (\ref{caso-inf}) is connected with the concentration
phenomenon (cf. \cite{le3} and the references therein).

For the corresponding analysis we start by combining (\ref{caso-inf}) with
(\ref{ija3})
\[
I(t)\geq ct\left(  \log\frac{1}{t}\right)  ^{1/2},\;t\in(0,\frac{1}{2}],
\]
to obtain%
\[
f^{\ast\ast}(t)-f^{\ast}(t)\preceq\frac{\left\|  f\right\|  _{Lip}}{\left(
\log\frac{1}{t}\right)  ^{1/2}},\text{ }t\in(0,\frac{1}{2}].
\]
Therefore, for $t\in(0,\frac{1}{2}],$ we have%
\begin{align*}
f^{\ast\ast}(t)-f^{\ast\ast}(1/2)  &  =\int_{t}^{1/2}\left(  f^{\ast\ast
}(s)-f^{\ast}(s)\right)  \frac{ds}{s}\\
&  \preceq\left\|  \left|  \nabla f\right|  \right\|  _{\infty}\int_{t}%
^{1/2}\frac{1}{\left(  \log\frac{1}{s}\right)  ^{1/2}}\frac{ds}{s}\\
&  \leq2\left\|  \left|  \nabla f\right|  \right\|  _{\infty}\left(  \log
\frac{1}{t}\right)  ^{1/2}.
\end{align*}
Thus, if $\lambda\left\|  \left|  \nabla f\right|  \right\|  _{\infty}%
^{2}\prec1,$%
\begin{align*}
\int_{0}^{1/2}e^{\lambda\left(  f^{\ast\ast}(t)-f^{\ast\ast}(1/2)\right)
^{2}}dt  &  \preceq\int_{0}^{1/2}e^{\left(  \log\frac{1}{t^{\lambda\left\|
\left|  \nabla f\right|  \right\|  _{\infty}^{2}}}\right)  }dt\\
&  =\int_{0}^{1/2}\frac{1}{t^{\lambda\left\|  \left|  \nabla f\right|
\right\|  _{\infty}^{2}}}dt<\infty.
\end{align*}
Moreover, since $f^{\ast\ast}$ is decreasing we have%
\begin{align*}
\int_{1/2}^{1}e^{\lambda\left(  f^{\ast\ast}(t)-f^{\ast\ast}(1/2)\right)
^{2}}dt  &  \leq\int_{1/2}^{1}e^{\lambda\left(  f^{\ast\ast}(1-t)-f^{\ast\ast
}(1/2)\right)  ^{2}}dt\\
&  =\int_{0}^{1/2}e^{\lambda\left(  f^{\ast\ast}(t)-f^{\ast\ast}(1/2)\right)
^{2}}dt.
\end{align*}
This readily implies the exponential integrability of $\left(  f(t)-f^{\ast
\ast}(1/2)\right)  :$%
\[
\int_{\mathbb{R}^{n}}e^{\lambda\left(  f(x)-f^{\ast\ast}(1/2)\right)  ^{2}%
}d\gamma_{n}(x)<\infty,
\]
and, in fact, we can readily compute the corresponding Orlicz norm.

In this fashion we are led to define a new space $L_{\log^{1/2}}(\infty
,\infty)(\mathbb{R}^{n},d\gamma_{n})$ by the condition\footnote{More
generally, the relevant spaces to measure exponential integrability to the
power $p$ are defined by%
\[
\sup(f^{\ast\ast}(t)-f^{\ast}(t))(\log\frac{1}{t})^{{1/p^{\prime}}}<\infty.
\]
}%
\begin{equation}
\left\|  f\right\|  _{L_{\log^{1/2}}(\infty,\infty)(\mathbb{R}^{n},d\gamma
_{n})}=\sup_{0<t<1}(f^{\ast\ast}(t)-f^{\ast}(t))\left(  \log\frac{1}%
{t}\right)  ^{1/2}<\infty. \label{definida}%
\end{equation}
Summarizing our discussion, we have%
\[
\left\|  f\right\|  _{L_{\log^{1/2}}(\infty,\infty)(\mathbb{R}^{n},d\gamma
_{n})}\preceq\left\|  \nabla f\right\|  _{L^{\infty}(\mathbb{R}^{n}%
,d\gamma_{n})}%
\]
and
\[
L_{\log^{1/2}}(\infty,\infty)(\mathbb{R}^{n},d\gamma_{n})\subset
e^{L^{2}(\mathbb{R}^{n},d\gamma_{n})}.
\]
The scale of spaces $\{L_{\log^{\alpha}}(\infty,\infty)\}_{\alpha\in R_{+}}$
is thus suitable to measure exponential integrability. When $\alpha=0$ we get
the celebrated $L(\infty,\infty)$ spaces introduced in \cite{bds}, which
characterize the rearrangement invariant hull of $BMO.$ The corresponding
underlying rearrangement inequality in the Euclidean case is the following
version of the John-Nirenberg lemma%
\[
f^{\ast\ast}(t)-f^{\ast}(t)\preceq\left(  f^{\#}\right)  ^{\ast}(t)
\]
where $f^{\#}$ is the sharp maximal operator used in the definition of $BMO$
(cf. \cite{bds} and \cite{km}).

In fact, in our context the $L(\infty,\infty)$ space is connected to the
exponential inequalities by Bobkov-G\"{o}tze \cite{bogo}. Proceeding as before
we see that (compare with \cite{bogo})%
\[
\left(  f^{\ast\ast}(t)-f^{\ast}(t)\right)  \preceq\left|  \nabla f\right|
^{\ast\ast}(t)\left(  \log\frac{1}{t}\right)  ^{-1/2},\;0<t<\frac{1}{2},
\]
from where if follows readily that $\left|  \nabla f\right|  \in e^{L^{2}%
}\Longrightarrow f\in L(\infty,\infty),$ and therefore if, moreover $\int
f=0,$ we can also conclude that $f\in e^{L}.$

\section{Symmetrization by truncation of entropy inequalities}

In this brief section we wish to indicate, somewhat informally, how our
methods can be extended to far more general setting. Let $(\Omega,\mu)$ be a
probability measure space. As in the literature, we consider the entropy
functional defined, on positive measurable functions, by%
\[
Ent(g)=\int g\log gd\mu-\int gd\mu\log\int gd\mu.
\]
Suppose for example that $Ent$ satisfies a lS inequality of order $1$ on a
suitable class of functions,%
\begin{equation}
Ent(g)\leq c\int\Gamma(g)d\mu. \label{abe}%
\end{equation}
Here $\Gamma$ is to be thought as an abstract gradient. We will make an
assumption that is not made in the literature but is crucial for our method to
work: we will assume that $\Gamma$ is *truncation friendly*, in the sense that
for any truncation of $f$ (see section \ref{sec3}) we have%
\begin{equation}
\left|  \Gamma(f_{h_{1}}^{h_{2}})\right|  =\left|  \Gamma(f)\right|
\chi_{\left\{  h_{1}<\left|  f\right|  \leq h_{2}\right\}  }. \label{abec}%
\end{equation}
While this is a non standard assumption, as we know, the usual gradients are
indeed *truncation friendly*. In order to continue we need the following
elementary result that comes from \cite{boze} (Lemma 2.2)
\begin{equation}
Ent(g)\geq-\log\left\|  g\right\|  _{0}\int gd\mu\label{abecd}%
\end{equation}
here $\left\|  g\right\|  _{0}=\mu\{g\neq0\}.$ Combining (\ref{abe}),
(\ref{abec}), (\ref{abecd}) it follows that%
\[
-\log\left\|  f_{h_{1}}^{h_{2}}\right\|  _{0}\int f_{h_{1}}^{h_{2}}d\mu\leq
c\int\left|  \Gamma(f)\right|  \chi_{\left\{  h_{1}<\left|  f\right|  \leq
h_{2}\right\}  }d\mu
\]%
\[
-\log\lambda_{f}(h_{1})\mu\{h_{1}<\left|  f(x)\right|  \leq h_{2}\}\leq
c\int\left|  \Gamma(f)\right|  \chi_{\left\{  h_{1}<\left|  f\right|  \leq
h_{2}\right\}  }d\mu
\]%
\[
\left(  -\log\lambda_{f}(h_{2})\right)  (h_{2}-h_{1})\lambda_{f}(h_{2})\leq
c\int_{\left\{  h_{1}<\left|  f\right|  \leq h_{2}\right\}  }\left|
\Gamma(f)\right|  d\mu
\]

Pick $h_{1}=f^{\ast}(s+h),$ $h_{2}=f^{\ast}(s),$ then%
\[
s(\log\frac{1}{s})\left(  f^{\ast}(s)-f^{\ast}(s+h)\right)  \leq
c\int_{\{f^{\ast}(s+h)<\left|  f\right|  \leq f^{\ast}(s)\}}\left|
\Gamma(f)\right|  d\mu.
\]
Thus,%
\[
s(\log\frac{1}{s})\frac{\left(  f^{\ast}(s)-f^{\ast}(s+h)\right)  }{h}%
\leq\frac{c}{h}\int_{\{f^{\ast}(s+h)<\left|  f\right|  \leq f^{\ast}%
(s)\}}\left|  \Gamma(f)\right|  d\mu.
\]
Therefore, following the analysis of Section \ref{secprinciple}, we find that,
for any Young's function $A,$ we have%
\[
A\left(  s(\log\frac{1}{s})(-f^{\ast})^{\prime}(s)\right)  \leq\frac{d}%
{ds}\left(  \int_{\left\{  \left|  f\right|  >f^{\ast}(s)\right\}  }A(\left|
\Gamma(f)\right|  )d\mu\right)  .
\]
Integrating, and using the Hardy-Littlewood-P\'{o}lya principle exactly as in
section \ref{secprinciple}, we obtain the following abstract version of the
P\'{o}lya-Szeg\"{o} principle%
\[
\int_{0}^{t}\left(  s(\log\frac{1}{s})(-f^{\ast})^{\prime}(s)\right)  ^{\ast
}(r)dr\leq\int_{0}^{t}\left|  \Gamma(f)\right|  ^{\ast}(r)dr.
\]
This analysis establishes a connection between entropy inequalities and
logarithmic Sobolev inequalities via symmetrization. In particular, our
inequalities extend the classical results to the setting of rearrangement
invariant spaces. For more details see \cite{mami}.

\end{document}